\documentclass{article}%
\usepackage{amssymb}
\usepackage{amsmath}
\usepackage{amsfonts}
\usepackage{graphicx}%
\newtheorem{theorem}{Theorem}

\begin{document}

\title{An Improved Energy-Area Inequality for Harmonic Maps Using Image Curvature}
\author{Simon P Morgan\\University of Minnesota}
\maketitle

\begin{abstract}
An ODE variational calculation shows that an image principle curvature ratio
factor, $\left(  \sqrt{\frac{\rho_{1}}{\rho_{2}}}+\sqrt{\frac{\rho_{2}}%
{\rho_{1}}}\right)  $, can raise the lower bound, 2(Image Area), on energy of
a harmonic map of a surface into $\mathbb{R}^{n}$. In certain situations,
including all radially symmetry harmonic maps, equality is achieved.

\end{abstract}

\section{Introduction}

The energy inequality for harmonic maps of surfaces into $\mathbb{R}^{n}$,
achieves equality if and only if the map is conformal and the image has zero
mean curvature [J][ES][EF][EL1][EL2]. We extend this connection between
principle curvatures and conformality of the harmonic map to non-minimal
surfaces. Note principle curvatures in $\mathbb{R}^{n}$ can be defined using
the vector valued second fundamental form. As the ratio of principle
curvatures deviates from 1, the minimal surface case, the map becomes higher
energy and less conformal. This enables the lower bound on energy to be raised
based on a principle curvature ratio term.

The results hinge on an ODE variational calculation using second order
derivatives of maps that are approximated to second degree (4) using principle
curvature information. This requires well known regularity results, e.g.:[L]
(thm 2.1.11), about harmonic maps of surfaces\ into $\mathbb{R}^{n}$ to be
sure that (4) is defined almost everywhere on images of harmonic maps.

\section{Physical intuition}

We say that as the image deviates more from being a minimal surface, more
energy is required of the harmonic map. A physical interpretation is to
consider an elastic sheet with low curvature in one direction and high
curvature in the other. To maintain equilibrium, the tension in the sheet in
the low curvature direction must be much greater than in the other direction,
thus contributing more to energy.

\section{The Inequality}

\begin{theorem}
If h is a degree 1 harmonic map from a smooth compact surface into
$\mathbb{R}^{n}$, $\rho_{1}$ and $\rho_{2}$ are principle curvatures of the
image of h, then%

\begin{equation}
Energy\geq\underset{image}{%
{\displaystyle\iint}
}\left(  \sqrt{\frac{\rho_{1}}{\rho_{2}}}+\sqrt{\frac{\rho_{2}}{\rho_{1}}%
}\right)  didj\geq2(area\text{ }of\text{ }image)
\end{equation}

whenever the integral makes sense on the image, taking 0/0=1. Also for a%
$>$%
0, a/0 =$\infty.$

Furthermore when the pull back to the domain of the directions of principle
curvatures are defined, we can say that equality is achieved on the left hand
side if and only if the pull back of the directions of principle curvatures
are orthogonal in the domain. This occurs for the radially symmetric case.
\end{theorem}

Remark: When $\rho_{1}$=$\rho_{2}$, the right hand inequality becomes
equality. The left hand inequality depends upon a version of the conformality
of the map. When maps are conformal and the image is a minimal surface, the
energy is twice the area of the image.

\textbf{Proof of Theorem 1}

Consider a local region of the image of a degree 1 harmonic map $h$ in
$\mathbb{R}^{3}$, where it has principle curvatures $\rho_{1}$ and $\rho_{2} $
This can be generalized to $\mathbb{R}^{n}$ using the vector valued second
fundamental form, so we shall continue the discussion only in $\mathbb{R}^{3}
$. Take intrinsic coordinates in the image $i$ and $j$ which are locally
orthonormal and parallel to the directions of principal curvatures, and having
origin at $h(p)$. Their pullbacks in the domain are directions $r$ and $s$,
having origin at $\mathit{p}$. Following the set up in figure 1, define local
orthonormal coordinates $\mathit{u},\mathit{v},$ on the domain. Finally the
range has orthonormal coordinates $\mathit{X}$\textit{, }and\textit{\ }%
$\mathit{Y}$\textit{,} tangential to\textit{\ }$\mathit{i}$\textit{\ }%
and\textit{\ }$\mathit{j}$\textit{\ }at\textit{\ }$\mathit{h(p)}$%
\textit{\ }and\textit{\ }$\mathit{Z}$\textit{, }parallel to principal
curvature radii at\textit{\ }$\mathit{h(p)}$. See figure 1.%
\begin{center}
\includegraphics[
width=3.34in
]%
{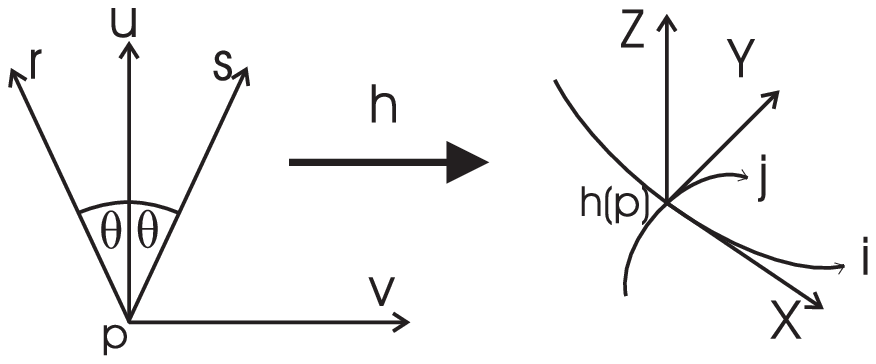}
\end{center}

\begin{center}
\textbf{Figure 1: Coordinate systems on domain and range}
\end{center}

Let the choice of \textit{u}, \textit{v}, \textit{s} and \textit{t} be such
that with linearization.

\begin{center}
$r=ucos(\theta)-v\sin(\theta)$

$s=ucos(\theta)+v\sin(\theta)$

$i_{r}(0,0)=X_{r}(0,0)=a\geq0,$ $j_{r}(0,0)=Y_{r}(0,0)=0$

$i_{s}(0,0)=Xs(0,0)=0,$ $j_{s}(0,0)=Y_{s}(0,0)=b\geq0$

$Z_{r}(0,0)=0,$ $Zs(0,0)=0,$

$\left[
\begin{array}
[c]{c}%
i\\
j
\end{array}
\right]  =\left[
\begin{array}
[c]{cc}%
a & 0\\
0 & b
\end{array}
\right]  \left[
\begin{array}
[c]{c}%
r\\
s
\end{array}
\right]  $

$\left[
\begin{array}
[c]{c}%
r\\
s
\end{array}
\right]  =\left[
\begin{array}
[c]{cc}%
\cos\left(  \theta\right)  & -\sin\left(  \theta\right) \\
\cos\left(  \theta\right)  & \sin\left(  \theta\right)
\end{array}
\right]  \left[
\begin{array}
[c]{c}%
u\\
v
\end{array}
\right]  $

$X_{u}^{2}(0,0)+X_{v}^{2}(0,0)=X_{r}^{2}(0,0)=i_{r}^{2}(0,0)$ and

$Y_{u}^{2}(0,0)+Y_{v}^{2}(0,0)=Y_{s}^{2}(0,0)=j_{s}^{2}(0,0)$
\end{center}

Now the energy of the map is%

\begin{equation}
\underset{domain}{%
{\displaystyle\iint}
}X_{u}^{2}+X_{v}^{2}+Y_{u}^{2}+Y_{v}^{2}+Z_{u}^{2}+Z_{v}^{2}dudv
\end{equation}

\begin{center}
=\bigskip$\underset{image}{%
{\displaystyle\iint}
}\frac{1}{\left\vert i_{r}j_{s}\sin\left(  2\theta\right)  \right\vert
}\left(  i_{r}^{2}+j_{s}^{2}\right)  didj$

=\bigskip$\underset{image}{%
{\displaystyle\iint}
}\frac{1}{\left\vert ab\sin\left(  2\theta\right)  \right\vert }\left(
a^{2}+b^{2}\right)  didj$
\end{center}

\begin{equation}
=\bigskip\underset{image}{%
{\displaystyle\iint}
}\frac{1}{\left\vert \sin\left(  2\theta\right)  \right\vert }\left(  \frac
{a}{b}+\frac{b}{a}\right)  didj
\end{equation}

We will now find the ratio $a/b$ in terms of the ratio of principal curvatures
on the image using a variational calculation. Let $X$, $Y$, and $Z $ be
orthonormal coordinates on the range as in figure 1. So we can write formulae
for the image coordinates, $X$, $Y$ and $Z$ in terms of image coordinates $i$
and $j$ and domain coordinates $r$ and $s$. The principle curvatures are
$\rho_{1}$ and $\rho_{2}$. To first degree $i=ar$ and $j=bs$ so $\rho_{1}ar$
and $\rho_{2}bs$ are linear approximations to the angle of rotation of the
image tangent plane, in directions of principle curvatures, with respect to
the tangent plane at the origin. This gives us the $Z$ coordinate, positive in
the direction of the radius of curvature $\rho_{1}$, based on a second degree approximation:%

\begin{equation}
Z(r,s)=\frac{1}{\rho_{1}}\left(  1-\cos\left(  \rho_{1}ar\right)  \right)
-\frac{1}{\rho_{2}}\left(  1-\cos\left(  \rho_{2}bs\right)  \right)  +o\left(
r^{3},r^{2}s,rs^{2},s^{3}\right)
\end{equation}

Now we place a smooth deformation field $\tau(X,Y)$ with compact support in
the $Z$ direction, and let it act with constant velocity, with time $t$. Now
we calculate the first variation of energy under the deformation.

\begin{center}
$\frac{d}{dt}Energy=\frac{d}{dt}\left(  \underset{image}{%
{\displaystyle\iint}
}\widetilde{Z}_{u}^{2}+\widetilde{Z}_{v}^{2}dudv\right)  =\underset
{}{\underset{image}{%
{\displaystyle\iint}
}\frac{1}{\left\vert \sin\left(  2\theta\right)  \right\vert }}\frac{d}%
{dt}\left(  \widetilde{Z}_{r}^{2}+\widetilde{Z}_{s}^{2}\right)  drds$

$\widetilde{Z}(u,v,t)=Z(u,v)+t\tau\left(  X(u,v),Y(u,v)\right)  $

$\widetilde{Z}_{u}(u,v,t)=Z_{u}(u,v)+t\tau_{u}\left(  X(u,v),Y(u,v)\right)  $

$\widetilde{Z}_{u}^{2}(u,v,t)=Z_{u}^{2}(u,v)+2tZ_{u}(u,v)\tau_{u}\left(
X(u,v),Y(u,v)\right)  +t^{2}\tau_{u}^{2}\left(  X(u,v),Y(u,v)\right)  $

$\frac{d}{dt}\left(  \widetilde{Z}_{u}^{2}(u,v,t)\right)  =2Z_{u}(u,v)\tau
_{u}\left(  X(u,v),Y(u,v)\right)  =2Z_{u}\tau_{u}$
\end{center}

Now we can write down the first variation of energy and simplify as follows:

\begin{center}
$\frac{d}{dt}Energy=2%
{\displaystyle\iint}
Z_{u}\tau_{u}+Z_{v}\tau_{v}dudv$
\end{center}

integrating by parts

\begin{center}
$=-2%
{\displaystyle\iint}
\left(  Z_{uu}+Z_{vv}\right)  \tau dudv$
\end{center}

Note that $\tau$ is smooth and bounded with compact support. This yields the
standard condition for a harmonic map that the coordinate functions are harmonic:

\begin{center}%
\begin{equation}
Z_{uu}+Z_{vv}=0
\end{equation}

\end{center}

Now using $Z_{rs}=o(r,s)$ and $Z_{s}=o(r^{2},rs,s^{2}),$ $Z_{r}=o(r^{2}%
,rs,s^{2})$ and from differentiating (4), we can obtain:

\begin{center}
$Z_{u}=r_{u}Z_{r}+s_{u}Z_{s}+o(r^{2},rs,s^{2})$

$Z_{uu}=\frac{\partial}{\partial u}\left(  r_{u}Z_{r}+s_{u}Z_{s}\right)
+o(r,s)$

$=r_{uu}Z_{r}+r_{u}Z_{ru}+s_{uu}Z_{s}+s_{u}Z_{su}+o(r,s)$

$Z_{uu}(0,0)=(r_{u}Z_{ru}+s_{u}Z_{su})$

$Z_{ru}=r_{u}Z_{rr}+s_{u}Z_{rs}+o(r,s)$

$Z_{ru}(0,0)=r_{u}Z_{rr}$
\end{center}

Now we can ignore higher order terms, giving:

\begin{center}%
\begin{equation}
Z_{uu}=r_{u}^{2}Z_{rr}+s_{u}^{2}Z_{ss}%
\end{equation}

\begin{equation}
Z_{vv}=r_{v}^{2}Z_{rr}+s_{v}^{2}Z_{ss}%
\end{equation}

\end{center}

Now we are evaluating $r_{u}$, $r_{v}$, $s_{u}$. and $s_{v}$ at $(0,0)$
Therefore we can use the linearization, as derivatives of quadratic terms on
the Taylor expansions of $r$ and $s$ will be zero at $(0,0)$. Also note that
in the Taylor expansion $\theta$ is a constant. So we can use:

\begin{center}
$r_{u}=\cos\left(  \theta\right)  +o(r,s),\qquad-r_{v}=\sin\left(
\theta\right)  +o(r,s),$

$s_{u}=\cos\left(  \theta\right)  +o(r,s),\qquad s_{v}=\sin\left(
\theta\right)  +o(r,s),$

$r_{u}^{2}(0,0)+r_{v}^{2}(0,0)=s_{u}^{2}(0,0)+s_{v}^{2}(0,0)$

$=cos^{2}\left(  \theta\right)  +\sin^{2}\left(  \theta\right)  =1$
\end{center}

substituting into (6) and (7) gives:

\begin{center}%
\[
Z_{uu}(0,0)+Z_{vv}(0,0)=Z_{rr}(0,0)+Z_{ss}(0,0)
\]

\end{center}

This point-wise calculation was for an arbitrary regular point, $(u_{0}%
,v_{0})$, where we set the local origin $(0,0)$ in our local coordinates. We
can repeat this for all regular points of the map. That is almost everywhere.
We can now say for the map to be energy stationary we obtain the condition in
terms of r and s:%

\begin{equation}
Z_{rr}+Z_{ss}=0
\end{equation}

Applying (4) to (8) gives:

\begin{center}
$Z(r,s)=\frac{1}{\rho_{1}}\left(  1-\cos\left(  \rho_{1}ar\right)  \right)
-\frac{1}{\rho_{2}}\left(  1-\cos\left(  \rho_{2}bs\right)  \right)  +o\left(
r^{3},r^{2}s,rs^{2},s^{3}\right)  $

$Z_{r}=\frac{\rho_{1}a}{\rho_{1}}\left(  \sin\left(  \rho_{1}ar\right)
\right)  +o(r^{2},rs,s^{2})$

$Z_{rr}=\rho_{1}a^{2}\left(  \cos\left(  \rho_{1}ar\right)  \right)  +o(r,s)$

$Z_{ss}=-\rho_{2}b^{2}\left(  \cos\left(  \rho_{2}bs\right)  \right)  +o(r,s)$

$\rho_{1}a^{2}\left(  \cos\left(  \rho_{1}ar\right)  \right)  -\rho_{2}%
b^{2}\left(  \cos\left(  \rho_{2}bs\right)  \right)  +o(r,s)=0$
\end{center}

Now taking $r$ and $s$ arbitrarily small by controlling the support of $\tau$
we can equate the constant terms resulting from the Taylor expansions to
obtain the relationship:%

\begin{equation}
\rho_{1}a^{2}-\rho_{2}b^{2}=0\Leftrightarrow\frac{a}{b}=\sqrt{\frac{\rho_{2}%
}{\rho_{1}}}%
\end{equation}

This gives us a lower bound on the energy in terms of the image and principal
curvatures on the image using (3) and (9), when the integrals make sense:%
\begin{align}
Energy  &  =\underset{image}{%
{\displaystyle\iint}
}\frac{1}{\left\vert \sin\left(  2\theta\right)  \right\vert }\left(
\sqrt{\frac{\rho_{1}}{\rho_{2}}}+\sqrt{\frac{\rho_{2}}{\rho_{1}}}\right)
didj\\
&  \geq\underset{image}{%
{\displaystyle\iint}
}\left(  \sqrt{\frac{\rho_{1}}{\rho_{2}}}+\sqrt{\frac{\rho_{2}}{\rho_{1}}%
}\right)  didj\\
&  \geq2(area\text{ }of\text{ }image)
\end{align}
$\blacksquare$

Now we can see that the condition $a=b$ corresponds to the principal
curvatures being equal, hence the surface being minimal. Also when $a=b$ and
the map is conformal, sin$\left(  2\theta\right)  \equiv1$ in (10). So we
conclude that energy equals twice the area of the image when the image is
minimal and the map is conformal as (10) and (11) become equalities.

Note the integral is defined when the quantity $\sqrt{\frac{\rho_{1}}{\rho
_{2}}}+\sqrt{\frac{\rho_{2}}{\rho_{1}}}$ is defined, i.e.: only for negatively
curved surfaces, not ruled or planar. Positive curvature does not arise in
images of harmonic maps into $\mathbb{R}^{n}$. In the planar case, $\rho
_{1}=\rho_{2}=0$ use $\sqrt{\frac{\rho_{1}}{\rho_{2}}}=1$. Note that this
corresponds to $a=b$ in (3) and (11) and (12) become equality.

In the radially symmetric case, $h:(r,\phi)\rightarrow(R,\Phi).$ $R$ and
$\Phi$ are only functions of one variable, $r$ and $\phi$ respectively.
Therefore in (10) $sin(2\theta)\equiv1$. Thus the inequality (11) becomes an
equality. This means that energy is completely determined by the image, its
area and curvatures.

\section{Acknowledgements}

Thanks to Robert Hardt, my thesis advisor, and to Michael Wolf and Robert
Gulliver for their insightful help.

\end{document}